\documentclass[10pt]{amsart}
\usepackage[margin=0.95in]{geometry}
\newtheorem{theorem}{Theorem}[section]

\usepackage{amsmath,amsthm,amsfonts,amssymb,amstext,euscript,verbatim,bbm}
\theoremstyle{definition}

\theoremstyle{remark}
\newtheorem{remark}[theorem]{Remark}

\numberwithin{equation}{section}

\begin{document}

\title[Gain of regularity for the relativistic collision operator]{Gain of regularity for the relativistic collision operator}

\author{Jin Woo Jang}
\address{Center for Geometry and Physics, Institute for Basic Science (IBS), Pohang 37673, Republic of Korea }
\email{jangjinw@ibs.re.kr}

\author{Seok-Bae Yun}
\address{Department of mathematics, Sungkyunkwan University, Suwon 440-746, Republic of Korea }
\email{sbyun01@skku.edu}




\begin{abstract}
We study a regularity property for the gain part of the relativistic Boltzmann collision operator. Our assumptions on the collisional scattering kernel cover the full range of generic hard and soft potentials with angular cut-off.
\end{abstract}
\allowdisplaybreaks
\maketitle
\section{introduction}
In this paper, we study the regularizing effect of the gain part of the collision operator of the relativistic Boltzmann equation:
\begin{eqnarray*}\label{RBE}
\partial_tf+\hat{p}\cdot\nabla_xf =Q(f,f).
\end{eqnarray*}
Here the momentum four vector $p^{\mu}$ is defined by $p^{\mu}=(p^0,p)$ where $p^0=\sqrt{1+|p|^2}$ and $p=(p_1,p_2,p_3)\in\mathbb{R}^3$ denote the energy and the momentum of a particle respectively. Here we have set the speed of light $c=1$ and the particle mass $m=1$. The normalized velocity  $\hat{p}$  is given by $\hat{p}= p/p^0=p/\sqrt{1+|p|^2}$. The velocity distribution function $f(x,p,t)$
represents the number density of the particle system on the phase point $(x,p)\in\mathbb{R}^3_x\times \mathbb{R}^3_p$ at time $t\in R_+$.
The collision operator $Q$ is given by
\begin{eqnarray*}
Q(f,h)(p)=Q^+(f,h)(p)-Q^-(f,h)(p).
\end{eqnarray*}
In the \textit{center-of-momentum} frame, the gain term $Q^+$ and the loss term $Q^-$ are written as \cite{DeGroot,Strain},
\begin{eqnarray*}
\displaystyle Q^+(f,h)=\int_{\mathbb{R}^3\times \mathbb{S}^2} v_{\o}\sigma(g,\theta)f(p^{\prime})h(q^{\prime})d\omega dq,\quad
\displaystyle Q^-(f,h)=\int_{\mathbb{R}^3\times \mathbb{S}^2} v_{\o}\sigma(g,\theta)f(p)h(q)d\omega dq,
\end{eqnarray*}
where $\sigma(g,\theta)$ is the scattering kernel, and the M{\o}ller velocity $v_{\o}$ is given by
\begin{eqnarray*}
v_{\o}=v_{\o}(p,q)=\sqrt{\left|\frac{p}{p_0}-\frac{q}{q_0}\right|^2-\left|\frac{p}{p_0}\times\frac{q}{q_0}\right|^2}
=\frac{g\sqrt{s}}{p_0q_0}.
\end{eqnarray*}
Here $s$ represents the square of the energy in the \textit{center-of-momentum} frame, $p+q=0$:
\begin{eqnarray*}
s=-(p^{\mu}+q^{\mu})(p_{\mu}+q_{\mu})
=2(p^0q^0-p\cdot q+1)\geq 0,
\end{eqnarray*}
and  $g$ denotes the relative momentum
\begin{eqnarray*}
g=\sqrt{(p^{\mu}-q^{\mu})(p_{\mu}-q_{\mu})}=\sqrt{2(p^0q^0-p\cdot q-1)}.
\end{eqnarray*}
Note that $s=g^2+4$.
The pre-collisional momentum pair $(p,q)$ and the post-collisional momentum pair $(p^{\prime},q^{\prime})$ are related by
\begin{equation}\label{p'}
\begin{split}
p^{\prime}=&\frac{p+q}{2}+\frac{g}{2}\left(\omega+(\gamma-1)(p+q)\frac{(p+q)\cdot\omega}{|p+q|^2}\right),\\
q^{\prime}=&\frac{p+q}{2}-\frac{g}{2}\left(\omega+(\gamma-1)(p+q)\frac{(p+q)\cdot\omega}{|p+q|^2}\right),
\end{split}\end{equation}
where $\gamma=(p^0+q^0)/\sqrt{s}$. The microscopic energy is given by
\begin{eqnarray*}
p^{\prime}_0=\frac{p_0+q_0}{2}-\frac{g}{2\sqrt{s}}\omega\cdot(p+q),\quad q^{\prime}_0=\frac{p_0+q_0}{2}-\frac{g}{2\sqrt{s}}\omega\cdot(p+q).
\end{eqnarray*}
The collision operator $Q$ satisfies
\begin{eqnarray}\label{symmetry}
\int_{\mathbb{R}^3}Q(g,g)(p)\left(1, p, p^0\right)dp=0,\quad \int_{\mathbb{R}^3}Q(g,g)(p)\ln f dp\leq 0,
\end{eqnarray}
which respectively leads to the conservation of mass, momentum and energy, and the celebrated H-theorem.


\subsection{Main results} We start with our main hypothesis on the collision kernel.\newline
\noindent{\bf $\bullet$ Hypothesis on $\sigma$}
The relativistic Boltzmann collision kernel $\sigma(g,\theta)$ is a non-negative function which depends only on the relative velocity $g$ and the scattering angle $\theta$. We assume that $\sigma$ satisfies
\begin{equation}
\label{hypothesis}\sigma(g,\theta)\lesssim g^a \sigma_0(\theta), \qquad(a>-3)
\end{equation}
where the angular kernel $\sigma_0(\theta)$ is nonnegative and bounded.

This hypothesis includes the full-range of generic hard- ($a\geq 0$) and soft-potential ($a< 0$) collisional kernels; for example, \eqref{hypothesis} with $a\geq 0$ includes the relativistic analogue of the hard-ball case with $a=1$ and the (elastic) neutrino gas with $a=2$, and \eqref{hypothesis} with $a<0$ includes the interactions of Israel particles with $a=-1$.
We now state our main result. In the following, $||\cdot||_{L^m_\nu}$ with $\nu\in \mathbb{R}$ denotes
$$||f||_{L^m_\nu(\mathbb{R}^3)}\equiv\left( \int_{\mathbb{R}^3} |f(p)|^m(p^0)^{\nu}dp\right)^{\frac{1}{m}},$$ where we abbreviate $f(x,p,t)$ by $f(p)$ for fixed $t$ and $x$.
\begin{theorem}\label{regularitymain}$[L^1\times L^2\rightarrow \dot{H^1}]$
$($\textit{Hard potentials}$)$  Assume that the scattering kernel $\sigma$ satisfies \eqref{hypothesis} with $a\geq 0$. Also, suppose that $f\in L^m_{\frac{m}{2}(a-1)}$ and $h\in L^n_{\frac{n}{2}(a-1)}$ with $\frac{1}{m}+\frac{1}{n}=\frac{3}{2}.$  Then the gain term $Q^+$ has the following regularizing property:
$$\|\nabla_p Q^{+}(f,h)\|_{L^2}\lesssim \|f\|_{L^m_{\frac{m}{2}(a-1)}}\|h\|_{L^n_{\frac{n}{2}(a-1)}}.$$
Note that, by choosing $m=1$ and $n=2$, we have
$$\|\nabla_p Q^{+}(f,h)\|_{L^2}\lesssim \|f\|_{L^1_{\frac{a-1}{2}}}\|h\|_{L^2_{a-1}}.$$
\end{theorem}

\begin{theorem} \label{regularitymain2}
$(1)$ $[L^2\times L^2\rightarrow \dot{H^1}]$ $($\textit{Hard potentials}$)$
Assume that the scattering kernel $\sigma$ satisfies \eqref{hypothesis} with $a\geq 0$. Suppose that $f,h\in L^2_{a+2+\varepsilon}$ for some $\varepsilon>0$.  Then the gain term $Q^+$ has the following regularizing property:
	$$\|\nabla_p Q^{+}(f,h)\|_{L^2}\lesssim \|f\|_{L^2_{a+2+\varepsilon}}\|h\|_{L^2_{a+2+\varepsilon}}.$$

\noindent $(2)$ $[L^m\times L^n\rightarrow \dot{H^1}]$ $($\textit{Soft potentials}$)$ Assume that the scattering kernel $\sigma$ satisfies \eqref{hypothesis} with $-3<a<0$. Suppose that $f\in L^{m}_{\left(\frac{|a|}{2}+1+\varepsilon\right)m}$ and
$h\in L^{n}_{\left(\frac{|a|}{2}+1+\varepsilon\right)n}$ with $\frac{1}{m}+\frac{1}{n}=1+\frac{a}{3}$
for some $\varepsilon>0$.  Then the gain term $Q^+$ has the following regularizing property:
\begin{align*}
\|\nabla_p Q^{+}(f,h)\|_{L^2}\lesssim \|f\|_{L^{m}_{\left(\frac{|a|}{2}+1+\varepsilon\right)m}}\|h\|_{L^{n}_{\left(\frac{|a|}{2}+1+\varepsilon\right)n}}.
\end{align*}
\end{theorem}

\begin{remark} The notation $\nabla_p$ above stands for the momentum derivatives with respect to $p\in\mathbb{R}^3$ only.\end{remark}\begin{remark}
Theorem \ref{regularitymain} holds only for the hard potential case with $a\geq0$, but it holds even if one of $f$ and $h$ lies in $L^1$ as in \cite{Andreasson,Lions,M-V,Wennberg2}.
In Theorem \ref{regularitymain2}, none of $f$ and $h$ can be a $L^1$ function, but the full-range of generic collision kernels with angular cut-off is covered.
\end{remark}

The smoothing estimate of the gain part of the collision operator was proved for the first time for regularized collision kernel by Lions in \cite{Lions},
who observed the Radon transform structure in the collision operator. 
Wennberg \cite{Wennberg2} then provided a proof which covers the physical collision kernel using the Carlemann representation of $Q^+$.
A much simplified proof using only the Fourier transform, at the cost of confining the function space into $L^2$,  was obtained in \cite{B-D,Lu}.
Mouhot and Villani generalized the estimates of \cite{B-D,Lu,Wennberg2} to weighted Sobolev spaces in \cite{M-V}.
Various  restrictions on the collision kernel were further relaxed by Jiang in \cite{Jiang}. Lions' proof was  generalized to the relativistic case by Andreasson in \cite{Andreasson}. Wennberg obtained a unified proof of such regularizing effects
for both classical and relativistic cases in \cite{Wennberg1} by showing that both can be deduced from the smoothing estimates of the generalized Radon transform.

In this paper, we provide a simplified proof of the regularizing estimate of the $Q^+$ operator from the relativistic Boltzmann equation in the spirit of \cite{B-D,Lu}.
Unlike \cite{B-D,Lu}, however, our proof gives the regularizing effect for the operator even when one of the functions lies in $L^1$ in the case of the full-range of the hard potentials $a\geq 0$ (Theorem \ref{regularitymain}) and shows the smoothing estimate for the full-range of generic hard and soft potentials
(Theorem \ref{regularitymain2}).


\section{Proof of Theorem 1.1}
We initially record the following pointwise estimates (see \cite{GS}, Lemma 3.1):
\begin{equation}
\label{pointwise g}
\frac{|p-q|}{\sqrt{p^0q^0}}\leq g\leq |p-q| \ \text{and}\ g\leq 2\sqrt{p^0q^0}.
\end{equation} Note that their definition of $``g"$ in \cite{GS} is actually a half of our $``g"$. Then we also have $s=4+g^2\lesssim p^0q^0$.

By the pre-post collisional change of variables $(p,q)\mapsto (p',q')$, we have
\begin{align}\label{1}
\begin{split}
\int_{\mathbb{R}^3}Q^+(f,h)e^{-ik\cdot p}dp&=\int_{\mathbb{R}^6\times\mathbb{S}^2}v_{\o}\sigma(g,\theta)f(p^{\prime})h(q^{\prime})e^{-ik\cdot p}d\omega dpdq\cr
&=\int_{\mathbb{R}^6\times \mathbb{S}^2}v_{\o}\sigma(g,\theta)f(p)h(q)e^{-ik\cdot p^{\prime}}d\omega dp^{\prime}dq^{\prime}\cr
&\equiv\int_{\mathbb{R}^6}v_{\o}\sigma(g,\theta)f(p)h(q)I(p,q,k) dp^{\prime}dq^{\prime},
\end{split}
\end{align}
where we define
$$I(p,q,k)\equiv \int_{\mathbb{S}^2}
e^{-ik\cdot \left\{\frac{p+q}{2}+\frac{g}{2}\left(\omega+(\gamma-1)(p+q)\frac{(p+q)\cdot\omega}{|p+q|^2}\right)\right\}}
d\omega .$$
Then $I(p,q,k)$ further satisfies that
\begin{align}\label{2}
\begin{split}
I
&=e^{-\frac{1}{2}i(p+q)\cdot k}
\int_{\mathbb{S}^2}e^{-i\frac{g}{2}\left\{\omega\cdot k+(\gamma-1)(p+q)\cdot k\frac{(p+q)\cdot\omega}{|p+q|^2}\right\}}d\omega\cr
&=e^{-\frac{1}{2}i(p+q)\cdot k}
\int_{\mathbb{S}^2}e^{-i\frac{g}{2}\left\{k+(\gamma-1)(p+q)\cdot k\frac{(p+q)}{|p+q|^2}\right\}\cdot\omega}d\omega\cr
&\equiv e^{-\frac{1}{2}i(p+q)\cdot k} II(p,q,k),
\end{split}
\end{align}
where $$ II(p,q,k)\equiv\int_{\mathbb{S}^2}e^{-i\frac{g}{2}\left\{k+(\gamma-1)(p+q)\cdot k\frac{(p+q)}{|p+q|^2}\right\}\cdot\omega}d\omega.$$ 
We then define $$A\equiv \frac{g}{2}\left\{k+(\gamma-1)(p+q)\cdot k\frac{(p+q)}{|p+q|^2}\right\},$$
and compute using the polar coordinate as
\begin{align*}
II&=\int_{\mathbb{S}^2}e^{-iA\cdot\omega}d\omega=\int_0^{2\pi}\left( \int_{0}^{\pi}\ e^{-i|A|\cos\psi}\sin\psi d\psi\right)d\vartheta\\
&=2\pi\int_{-1}^{1}e^{-i|A|s}ds=\frac{2\pi}{i|A|}(e^{i|A|}-e^{-i|A|})\\
&=4\pi\frac{\sin\left(\left|\frac{g}{2}\left\{k+(\gamma-1)(p+q)\cdot k\frac{(p+q)}{|p+q|^2}\right\}\right|\right)}{\left|\frac{g}{2}\left\{k+(\gamma-1)(p+q)\cdot k\frac{(p+q)}{|p+q|^2}\right\}\right|}.
\end{align*}
Since $|\sin(|A|)|\leq 1$, this gives
\begin{align}\label{3}
|II|\leq\frac{8\pi}{g\left|k+(\gamma-1)(p+q)\cdot k\frac{(p+q)}{|p+q|^2}\right|}.
\end{align}
Combining (\ref{1}), (\ref{2}) and  (\ref{3}), we arrive at
\begin{align}\label{4}
\bigg|\int_{\mathbb{R}^3}Q^+(f,h)e^{-ik\cdot p}dp\bigg|\leq C\left|\int_{\mathbb{R}^6}v_{\o}\sigma(g,\theta)g^{-1}
\frac{e^{-\frac{1}{2}i k\cdot(p+q)}f(p)h(q) }{\left|k+(\gamma-1)(p+q)\cdot k\frac{(p+q)}{|p+q|^2}\right|}dpdq\right|.
\end{align}
We then observe from $\displaystyle|x|=\max_{|y|=1}|y^{\top}x|$ that
\begin{align*}
\left|k+(\gamma-1)(p+q)\cdot k\frac{(p+q)}{|p+q|^2}\right|&=\big|\left(I+(\gamma-1)A\otimes A\right)k\big|\cr
&=\max_{|y|=1}\big|y^{\top}\left(I+(\gamma-1)A\otimes A\right)k\big|\cr
&\geq\left|\left(\frac{k}{|k|}\right)^{\top}\left(I+(\gamma-1)A\otimes A\right)k\right|\cr
&=\frac{1}{|k|}\left(|k|^2+(\gamma-1)\left\{A\cdot k\right\}^2\right)\cr
&\geq |k|,
\end{align*}
where $A=\frac{p+q}{|p+q|}$. The last inequality holds because $$\gamma-1=\frac{p^0+q^0-\sqrt{s}}{\sqrt{s}}=\frac{|p+q|^2}{\sqrt{s}(p^0+q^0+\sqrt{s})}>0.$$
Therefore we obtain from (\ref{4}) that
\begin{align}\label{ineq1}
\left|\int_{\mathbb{R}^3}Q^+(f,h)e^{-ik\cdot p}dp\right|\leq\frac{C}{|k|}\bigg|\int_{\mathbb{R}^6}v_{\o}\sigma(g,\theta)g^{-1}e^{-\frac{1}{2}i k\cdot(p+q)}f(p)h(q) dpdq\bigg|.
\end{align}
Note that 
$g\lesssim \sqrt{p^0q^0}$ and $s\lesssim p^0q^0$. Since $\sigma_0$ is bounded, we then have
\begin{equation}\label{vphi1}
v_{\o}\sigma(g,\theta)g^{-1}\lesssim \frac{g\sqrt{s}}{p^0q^0}g^{a-1}\approx g^a\frac{\sqrt{s}}{p^0q^0}\lesssim (p^0q^0)^{\frac{a}{2}}\frac{1}{\sqrt{p^0q^0}},
\end{equation}
because of the pointwise estimates \eqref{pointwise g} and that $a\geq 0$.
Therefore,
\begin{align*}\
\displaystyle|k|\left|\widehat{Q^{+}(f,h)}(k)\right|\lesssim\left|\int_{\mathbb{R}^6}e^{-\frac{1}{2}ik\cdot(p+q)}(p^0q^0)^{\frac{a-1}{2}}f(p)h(q)dpdq\right|\approx\left|\widehat{f_a}\left(\frac{k}{2}\right)\right|\left|\widehat{h_a}\left(\frac{k}{2}\right)\right|,\notag
\end{align*}
where $f_a(p)\equiv (p^0)^{\frac{a-1}{2}}f(p)$ and $h_a(q)\equiv (q^0)^{\frac{a-1}{2}}h(q)$. On the other hand, Plancherel's identity gives
$$
||\nabla_p Q^+||^2_{L^2}=\left|\left|\widehat{\nabla_p Q^+}\right|\right|^2_{L^2}=\left|\left||k|\widehat{Q^+}\right|\right|^2_{L^2}.
$$
Altogether, we obtain that
\begin{align*}
||\nabla_p Q^+(f,h)||^2_{L^2}
\lesssim\int_{\mathbb{R}^3} \left|\widehat{f_a}\left(\frac{k}{2}\right)\widehat{h_a}\left(\frac{k}{2}\right)\right|^2dk
\approx \int_{\mathbb{R}^3}\left|\widehat{f_a*h_a}\left(\frac{k}{2}\right)\right|^2dk
\approx \|\widehat{f_a*h_a}\|^2_{L^2}\approx\|f_a*h_a\|^2_{L^2},
\end{align*}
where we use the convolution formula for the Fourier transform and Plancherel's identity. Then, by Young's convolution inequality, we finally have
$$
||\nabla_p Q^+(f,h)||^2_{L^2}\lesssim \|f_a\|^2_{L^m}\|h_a\|^2_{L^n},\quad\quad \left(\frac{1}{m}+\frac{1}{n}=\frac{3}{2}\right).
$$
That is,
$$
||\nabla_p Q^+(f,h)||_{L^2}\lesssim \|f\|_{L^m_{\frac{m}{2}(a-1)}}\|h\|_{L^n_{\frac{n}{2}(a-1)}}.$$
This completes the proof for Theorem \ref{regularitymain}.

\section{Proof of Theorem 1.2} \subsection{Hard potential case  $(a\geq 0)$:}
 We start with \eqref{ineq1} and \eqref{vphi1}.
By taking the $L^2$ norm of the left-hand side of \eqref{ineq1} and using Plancherel's identity and \eqref{vphi1}, we have
$$\|\nabla_p Q^{+}(f,h)\|^2_{L^2}\lesssim \int_{\mathbb{R}^3}\left|\int_{\mathbb{R}^6}e^{-\frac{1}{2}ik\cdot(p+q)}(p^0q^0)^{\frac
{a-1}{2}}f(p)h(q)dpdq\right|^2dk.$$
Then we take the following change of variables:
\begin{align}\label{CV}
x=\frac{p+q}{2},\quad y=p-q,
\end{align}
and denote
\[
H_y(x)=f\left(x+\frac{y}{2}\right)h\left(x-\frac{y}{2}\right)\left(1+\left|x+\frac{y}{2}\right|^2\right)^{\frac
{a-1}{4}}\left(1+\left|x-\frac{y}{2}\right|^2\right)^{\frac
{a-1}{4}}
\]
to derive
\begin{align*}
\|\nabla_p Q^{+}\|^2_{L^2}\lesssim
\int_{\mathbb{R}^3}\left|\int_{\mathbb{R}^3}\left\{\int_{\mathbb{R}^3}e^{-ik\cdot x}H_y(x)dx\right\}dy\right|^2dk=\int_{\mathbb{R}^3}\left|\int_{\mathbb{R}^3}\widehat{H_y}(k)dy\right|^2dk.
\end{align*}
Then by Cauchy-Schwarz inequality,
\begin{align*}
\|\nabla_p Q^{+}\|^2_{L^2}&\lesssim\int_{\mathbb{R}^3}\left(\int_{\mathbb{R}^3}|\widehat{H_y}(k)|(1+|y|^2)^{\frac{3+\varepsilon}{4}}(1+|y|^2)^{-\frac{3+\varepsilon}{4}}dy\right)^2dk\cr
&\lesssim\int_{\mathbb{R}^3}\left\{\int_{\mathbb{R}^3}|\widehat{H_y}(k)|^2(1+|y|^2)^{\frac{3+\varepsilon}{2}}dy\int_{\mathbb{R}^3}(1+|y|^2)^{-\frac{3+\varepsilon}{2}}dy\right\}dk\cr
&=C_{\varepsilon}\int_{\mathbb{R}^3}\left\{\int_{\mathbb{R}^3}|\widehat{H_y}(k)|^2(1+|y|^2)^{\frac{3+\varepsilon}{2}}dy\right\}dk,
\end{align*}
and by Fubini's theorem and Plancherel's identity, we have
\begin{align*}
\|\nabla_p Q^{+}\|^2_{L^2}
&\lesssim \int_{\mathbb{R}^3}\left\{\int_{\mathbb{R}^3}|H_y(x)|^2dx\right\}(1+|y|^2)^{\frac{3+\varepsilon}{2}}dy\cr
&\approx\int_{\mathbb{R}^6}(p^0q^0)^{
a-1}|f(p)|^2|h(q)|^2(1+|p-q|^2)^{\frac{3+\varepsilon}{2}}dpdq\cr
&\lesssim \int_{\mathbb{R}^6}(p^0q^0)^{a-1+3+\varepsilon}|f(p)|^2|h(q)|^2dpdq\cr
&\lesssim\|f\|^2_{L^2_{a+2+\varepsilon}}\|h\|^2_{L^2_{a+2+\varepsilon}},
\end{align*}
where we used
\begin{align}\label{p-q}
1+|p-q|^2\lesssim (p^0q^0)^2.
\end{align}
This completes the proof for the hard-potential case with $a\geq 0$.

\subsection{Soft potential case $(a<0)$:} Again, we start from \eqref{ineq1}.
In this case, the scattering kernel becomes more singular in $g$, and we use the following sharp inequality instead of \eqref{vphi1}:
	$$v_{\o}\sigma(g,\theta)g^{-1}\lesssim \frac{g\sqrt{s}}{p^0q^0} g^{a-1}\lesssim \left(\frac{|p-q|}{\sqrt{p^0q^0}}\right)^{a}\frac{1}{\sqrt{p^0q^0}},$$ where we used $v_{\o}=\frac{g\sqrt{s}}{p^0q^0}$, $a<0$, $g\gtrsim \frac{|p-q|}{\sqrt{p^0q^0}}$, and $s\lesssim p^0q^0$ from the pointwise estimates \eqref{pointwise g}. By taking the $L^2$ norm of the left-hand side of \eqref{ineq1} and using Plancherel's identity, we have
	$$\|\nabla_p Q^{+}\|^2_{L^2}\lesssim \int_{\mathbb{R}^3}\left|\int_{\mathbb{R}^6}e^{-\frac{1}{2}ik\cdot(p+q)}\left(\frac{|p-q|}{\sqrt{p^0q^0}}\right)^{a}\frac{1}{\sqrt{p^0q^0}}f(p)h(q)dpdq\right|^2dk.$$
	Then we take the same change of variables as in (\ref{CV})
	and define
	\[
	S_y(x)=|y|^{a}f\left(x+\frac{y}{2}\right)h\left(x-\frac{y}{2}\right)\left\{\left(1+\left|x+\frac{y}{2}\right|^2\right)\left(1+\left|x-\frac{y}{2}\right|^2\right)\right\}^{-\frac{a+1}{4}},
	\]
to obtain
	\begin{align*}
	\|\nabla_p Q^{+}\|^2_{L^2}
	=\int_{\mathbb{R}^3}\left|\int_{\mathbb{R}^3}\left\{\int_{\mathbb{R}^3}e^{-ik\cdot x}S_y(x)dx\right\}dy\right|^2dk
	=\int_{\mathbb{R}^3}\left|\int_{\mathbb{R}^3}\widehat{S_y}(k)dy\right|^2dk.
	\end{align*}
	Then by Cauchy-Schwarz inequality, we have
	\begin{align*}
	\|\nabla_p Q^{+}\|^2_{L^2}&\lesssim\int_{\mathbb{R}^3}\left(\int_{\mathbb{R}^3}|\widehat{S_y}(k)|(1+|y|^2)^{\frac{3+\delta}{4}}(1+|y|^2)^{-\frac{3+\delta}{4}}dy\right)^2dk\cr
	&\lesssim\int_{\mathbb{R}^3}\left\{\int_{\mathbb{R}^3}|\widehat{S_y}(k)|^2(1+|y|^2)^{\frac{3+\delta}{2}}dy\int_{\mathbb{R}^3}(1+|y|^2)^{-\frac{3+\delta}{2}}dy\right\}dk\cr
	&=C_{\delta}\int_{\mathbb{R}^3}\left\{\int_{\mathbb{R}^3}|\widehat{S_y}(k)|^2(1+|y|^2)^{\frac{3+\delta}{2}}dy\right\}dk.
	\end{align*}

Applying Fubini's theorem and Plancherel's identity, we obtain
	\begin{align*}
	\|\nabla_p Q^{+}(f,h)\|^2_{L^2}
	&\lesssim \int_{\mathbb{R}^3}\left\{\int_{\mathbb{R}^3}|S_y(x)|^2dx\right\}(1+|y|^2)^{\frac{3+\delta}{2}}dy\cr
	&\approx\int_{\mathbb{R}^6}\left(\frac{|p-q|}{\sqrt{p^0q^0}}\right)^{2a}\frac{1}{p^0q^0}|f(p)|^2|h(q)|^2(1+|p-q|^2)^{\frac{3+\delta}{2}}dpdq\cr
&\lesssim\int_{\mathbb{R}^6}\frac{|f(p)|^2|h(q)|^2}{|p-q|^{2|a|}}(p^0q^0)^{|a|-1+3+\delta}dpdq\cr
&\approx\int_{\mathbb{R}^6}\frac{|f(p)|^2|h(q)|^2}{|p-q|^{2|a|}}(p^0q^0)^{|a|+2+\delta}dpdq.
\end{align*}
We then recall the Hardy-Littlewood-Sobolev inequality:
\begin{align*}
\|\nabla_p Q^{+}(f,h)\|^2_{L^2}&\lesssim \big\||f(p)|^2(p^0)^{|a|+2+\delta}\big\|_{L^{m'}}\big\||h(p)|^2(p^0)^{|a|+2+\delta}\big\|_{L^{n'}}
\quad \left(\frac{1}{m'}+\frac{1}{n'}=2-\frac{2|a|}{3}\right)\cr
&\approx\|f\|^2_{L^{2m'}_{\left(|a|+2+\delta\right)m'}}\|h\|^2_{L^{2n'}_{\left(|a|+2+\delta\right)n'}}.
\end{align*}
Finally, we rewrite $2m'=m$, $2n'=n$ and $\frac{\delta}{2}=\varepsilon$ to complete the proof.
\section{Acknowledgement}J. W. Jang was supported by
the Korean IBS project IBS-R003-D1. S.-B. Yun was supported by Basic Science Research Program through the National Research Foundation of Korea (NRF) funded by the Ministry of Education (NRF-2016R1D1A1B03935955).
\bibliographystyle{amsplain}

\end{document}